%%%%%%%%%%%%%%%%%%%%%%%%%%%%%%%%%%%%%%%%%%%%%%%%%%%%%%%%%%%%%%%%%%%%%%%%%%%%%%%
%       ZERO DIFFUSION-DISPERSION LIMITS FOR SCALAR CONSERVATION LAWS  
%       by Cezar Kondo and Philippe G. LeFloch 
%       SIAM Journal on Mathematical Analysis  
%%%%%%%%%%%%%%%%%%%%%%%%%%%%%%%%%%%%%%%%%%%%%%%%%%%%%%%%%%%%%%%%%%%%%%%%%%%%%% 
\input amstex
\documentstyle{amsppt}
\magnification 1100
\TagsOnRight
\NoBlackBoxes
\nologo  
\vsize=8.8 truein
\NoBlackBoxes 
%----------------------------------------------------------------------

\newif\iftitle

\def\leaderfill{ \leaders\hbox to1em{\hss.\hss}\hfill}

\newwrite\ind
\immediate\openout\ind=indice
\immediate\write\ind{\string\input }
\immediate\write\ind{fmtind}

\def\beginsection#1#2#3\par{\vskip0pt plus.3\vsize\penalty10000
    \vskip0pt plus -.3\vsize
\ifcase#1
\bigskip\vskip\parskip %\leftline
\item{\bf #2}  \bf #3\smallskip\or
\noindent{\item{\bf#2}\it #3\hfil\break\vskip0pt}
\else %\leftline
{\item{\it#2} \it #3}\fi
\message{#2#3}
\edef\mywrite##1##2##3{%
  \noexpand\writeind {##1}{##2}{##3} }%
\mywrite{#1}{#2}{#3}%
\ifcase#1\mark{#2\quad#3}\fi
%\np
\rm}

\def\writeind #1#2#3 {
\write\ind {
 \string\llap{\hbox to30pt{#2\string\hfil}}
                 #3\string\leaderfill \string\quad
 {\string\oldnos \folio}\string\break }
 }

\def\label#1{\smallskip\par\noindent\hangindent30pt #1}

\def\ra{\rangle}
\def\la{\langle}

\def\label#1{\smallskip\par\noindent\hangindent30pt#1}
%====================================================================================================

\def\RR    {{R\!\!\!\!\!I~}} 
\def\Rd   {{\RR^d}} 
\def\del   {{\partial}}
\def\eps   {\epsilon} 

\def\ra    {\big\rangle}
\def\la    {\big\langle} 
%______________________ 
\define  \refBLP  {1}
\define  \refBL   {2}
\define  \refCL   {3}
\define  \refDiPerna {4} 
\define  \refHLone   {5} 
\define  \refHLtwo   {6} 
\define  \refHLthree {7} 
\define  \refJMS     {8} 
\define  \refKL      {9} 
\define  \refKruzkov {10}
\define \refLax      {11} 
\define  \refLL      {12} 
\define  \refLeFlochone {13} 
\define  \refLeFlochtwo {14} 
\define  \refLN         {15}
\define  \refLR         {16} 
\define  \refLT         {17} 
\define  \refMurat      {18} 
\define  \refSchonbek   {19} 
\define  \refSzepessy   {20} 
\define  \refTartar     {21}

%=======================================================================================================

%%%%%%%%%%%%%%%%%%%%%%%%%%%%%%%%%%%%%%%%%%%%%%%%%%%%%%%%%%%%%%%%%%%%%%
\topmatter 
\title 
ZERO~DIFFUSION-DISPERSION~LIMITS\\
FOR SCALAR~CONSERVATION LAWS
\endtitle
\thanks Published in SIAM Math. Anal. 33 (2002), 1320--1329. 
\endthanks
\endtopmatter
%-----------------------------------------------------------------------
\bigskip

{\baselineskip11pt 

\centerline{\bf Cezar I. Kondo}

\medskip
 
\centerline{Universidade Federal de S$\tilde{a}$o Carlos} 
\centerline{Departamento de Matem\'atica, Caixa Postal 676} 
\centerline{13565-905, S$\tilde{a}$o Carlos-SP, Brazil.}  
}
\bigskip

\centerline{and}

\bigskip

\centerline{\bf
Philippe G. LeFloch\footnote"$^*$"{Former address : Ecole Polytechnique, Palaiseau, France (until 2004). 
\newline 
AMS {\it Mathematics Classification (2000)}\/. Primary: 35L65. Secondary: 76N10. 
\newline This work was done when the first author was visiting  
the Ecole Polytechnique during the Academic Year 1999--2000
with a financial support from PICD-Capes (Brazil). % and the European research project Alpha.
}
}

\medskip

\centerline{Laboratoire Jacques-Louis Lions}
\centerline{Centre National de la Recherche Scientifique}
\centerline{Universit\'e de Paris 6}
\centerline{4, Place Jussieu, 75252 Paris, France.}
\centerline{Email: LeFloch\@ann.jussieu.fr} 

\bigskip

\medskip

\bigskip\bigskip{\narrower\smallskip\noindent
\baselineskip11pt\eightpoint{\bf Abstract.\/} 
We consider solutions of 
hyperbolic conservation laws regularized with vanishing diffusion and dispersion terms. 
Following a pioneering work by Schonbek, we establish the convergence of the regularized 
solutions toward discontinuous solutions of the hyperbolic conservation law. 
The proof relies on the method of compensated compactness in the $L^2$ setting. 
Our result improves upon Schonbek's earlier results and provides an optimal condition on 
the balance between the relative sizes of the diffusion and the dispersion parameters. 
A convergence result is also established for multi-dimensional conservation laws  
by relying on DiPerna's uniqueness theorem for entropy measure-valued solutions. 

\medskip

\noindent{\it Key words and phrases:} conservation law, shock wave, 
entropy solution, measure-valued solution, diffusion, dispersion, 
singular limit, a priori estimate.  
\smallskip}

\vfill
\eject

%=======================================================================================================

\heading{1. Introduction}
\endheading

We study here the convergence of solutions of the partial differential equation
($\eps \to 0+$, $\delta= \delta(\eps) \to 0$) 
$$ 
u_t + f(u)_x 
= \eps \, u_{xx} + \delta \, u_{xxx}, \qquad u = u^\eps(x,t), \, x \in \RR, \, t \geq 0, 
\tag 1.1 
$$ 
toward weak solutions of the corresponding hyperbolic conservation laws: 
$$ 
u_t + f(u)_x = 0, \qquad u=u(x,t), \, x \in \RR, \, t \geq 0, 
\tag 1.2 
$$ 
where the flux $f: \RR \to \RR$ is a smooth function with (at most) 
linear growth at infinity, that is, for some $M > 0$ 
$$ 
|f'(u)| \leq M, \qquad u \in \RR. 
$$ 

Equations of the form (1.1)-(1.2) arise in fluid dynamics when
both viscosity (diffusion) and capillarity (dispersion) play a role. 
The diffusion $\eps$ smoothes  out the discontinuous solutions of (1.2), while 
the dispersion $\delta$ causes high-frequency oscillations. 

In this paper, we establish that the solutions $u^\eps$ of (1.1) 
converge toward a weak solution of (1.2) provided 
$$
\delta = O(\eps^2).
\tag 1.3  
$$
When the stronger condition 
$$
\delta = o(\eps^2) 
\tag 1.4
$$
holds, we prove that the limit coincides with the entropy solution 
determined by Kruzkov's theory \cite{\refKruzkov}. We point out that 
these conditions are sharp since, in the limiting case,  
$$
\delta = K \, \eps^2 \qquad \text { for some } K \in \RR, 
\tag 1.5 
$$
limiting solutions may violate Kruzkov's entropy conditions \cite{\refJMS, \refHLone, 
\refLeFlochone, \refBL}. 
Furthermore, when (1.3) is violated, the solutions are highly-oscillatory 
and fail to converge in any strong topology as noted by Lax and Levermore \cite{\refLL}. 
(See also Lax \cite{\refLax}. )

The singular limit problem above was first tackled by Schonbek \cite{\refSchonbek}, 
who established the optimal rate (1.3) for Burgers equation, that is,  
$$
f(u) = {u^2 \over 2}, 
$$
and for the class of flux-functions  
$$
f(u) = {u^{2p+1} \over 2p+1}, \qquad p \geq 1. 
$$ 
She also gave a convergence result for general fluxes with quadratic growth at infinity, 
however under the  stronger condition on $\delta= O(\eps^3)$. 
As another important contribution in \cite{\refSchonbek}, Schonbek introduces 
a generalization of the method of compensated compactness (Tartar \cite{\refTartar} 
and Murat \cite{\refMurat}) allowing to handle sequences that are bounded in 
$L^p$ for finite $p>1$ only.  Next, following \cite{\refSchonbek}, LeFloch 
and Natalini \cite{\refLN} 
studied equations like (1.1) but with nonlinear (even singular) diffusion, 
and established strong convergence results toward entropy solutions of (1.2). 
See also a convergence result for systems in Hayes and LeFloch \cite{\refHLtwo}.

In the second part of this paper, 
we also deal with the convergence of solutions of multi-dimensional equations 
similar to (1.1)-(1.2). For multi-dimensional equations, the compensated compactness 
method no longer applies and the proofs are based instead on DiPerna's uniqueness theory 
for entropy measure-valued solutions (DiPerna \cite{\refDiPerna}, Szepessy 
\cite{\refSzepessy}, and Kondo and LeFloch \cite{\refKL}). 
Our approach is  similar to Correia and LeFloch \cite{\refCL} where nonlinear 
diffusion terms are treated under a strong assumption on the ratio of the dispersion to the diffusion. 

To summarize, the main contribution in the present paper is the derivation of  
a priori estimates  (Theorems 2.1 and 3.1) which cover general flux-functions 
(with at most linear growth at infinity) and lead to an optimal condition on the balance 
between the diffusion and the dispersion (Theorems 2.2 and 3.2). 

Further material on classical and nonclassical entropy solutions 
generated by diffusive-dispersive limits can be found in 
\cite{\refBLP, \refBL, \refHLone, \refHLtwo, \refHLthree, \refJMS, 
\refLeFlochone, \refLeFlochtwo, \refLN, \refLR, \refLT, \refSchonbek}

%============================================================== 

\heading{2. One-Dimensional Conservation Laws} 
\endheading

Consider a family $u^\eps$ of smooth solutions to  
$$ 
u_t + f(u)_x = \eps \, u_{xx} + \delta \, u_{xxx},
\qquad u = u^\eps(x,t), 
\tag 2.1 
$$ 
$$
u(x,0)= u_0^\eps(x),  \qquad x \in \RR, 
\tag 2.2
$$
where $\eps \to 0+$ and $\delta = \delta(\eps) \to 0$. 
Under suitable conditions on the initial data $u_0^\eps: \RR \to \RR$, the solutions 
(and their derivatives) decay rapidly at infinity, so that all the a priori 
estimates given below are rigorously justified. 
We want to show that the solution of (2.1)-(2.2) converges toward a weak solution
of the problem 
$$ 
u_t + f(u)_x = 0, \qquad u = u^\eps(x,t), 
\tag 2.3 
$$ 
$$
u(x,0)= u_0(x),  \qquad x \in \RR, 
\tag 2.4
$$
where $u_0: \RR \to \RR$ is a given initial data. A minimum requirement is 
the weak convergence (for instance in $L^2(\RR)$) 
$$
u_0^\eps \rightharpoonup u_0, 
$$  
which is always assumed throughout this paper. 
The following convergence theorem covers both cases where the diffusion 
are in balance or dominates the dispersion. 

\proclaim{Theorem 2.1} Suppose that the flux-function the flux-function 
$f$ is Lipschitz continuous on $\RR$
and that the initial data $u_0$ belong to $L^2(\RR)$. 
Then the solution $u^\eps$ of $(2.1)$-$(2.2)$ satisfies the following a priori estimates: 
$$
\|u^\eps(t)\|_{L^2(\RR)} 
\leq \|u^\eps_0\|_{L^2(\RR)}, \qquad t \geq 0, 
\tag 2.5a
$$
$$
\sqrt { 2 \, \eps } \, \|u^\eps_x\|_{L^1\bigl(\RR^+,L^2(\RR)\bigr)}
\leq \|u^\eps_0\|_{L^2(\RR)}, 
\tag 2.5b 
$$ 
$$
\sqrt{\delta} \, \|u_x^\eps(t)\|_{L^2(\RR)} 
\leq \sqrt{ 2 \, \|f'\|_\infty} \, \|u^\eps_0\|_{L^2(\RR)} 
        + \sqrt{\delta} \,  \|u_{0 x}^\eps\|_{L^2(\RR)}, \qquad t \geq 0, 
\tag 2.5c
$$
and 
$$
\sqrt{ \eps \, \delta} \, \|u_{xx}\|_{L^1\bigl(\RR^+,L^2(\RR)\bigr)}
\leq \sqrt{ 2 \, \|f'\|_\infty} \, \|u^\eps_0\|_{L^2(\RR)} 
        + \sqrt{\delta} \,  \|u_{0 x}^\eps\|_{L^2(\RR)}. 
\tag 2.5d 
$$
\endproclaim 

%______________________________________________________________ 

\demo{Proof} Throughout the calculation and for simplicity, we omit the upper-index $\eps$. 
To any smooth function $U: \RR \to \RR$ we can associate a ``flux''  
$F:\RR \to \RR$ by $F'(u) = U'(u) f'(u)$, $u \in \RR$. 
Multiplying (2.1) by $U'(u)$ we find
$$
U(u)_t + F(u)_x 
= \eps \, (U'(u) u_x)_x - \eps \, U''(u) \, u^2_x
+ \delta \, \bigl(U'(u) \, u_{xx} \bigr)_x - \delta U''(u) \, u_x \, u_{xx}. 
$$ 
Integrating over the whole space, it follows that 
$$
\aligned 
{d \over dt} \int_\RR U(u) \, dx 
+ 
\eps \, \int_\RR U''(u) \, u^2_x \, dx 
& =
\delta \, \int_\RR U''(u) \, \bigg({u^2_x \over 2}\bigg)_x \, dx
\\
&= 
- {\delta \over 2} \, \int_\RR U'''(u) \, u^3_x \, dx. 
\endaligned 
\tag 2.6 
$$
Integrating in time over some interval $(0,t)$, we arrive at the general identity: 
$$
\int_\RR U(u(t)) \, dx + \eps \, \int_0^t \int_\RR U''(u) \, u^2_x \, dx dt
=  \int_\RR U(u_0) \, dx 
- {\delta \over 2} \, \int_0^t \int_\RR U'''(u) \, u^3_x \, dx dt.
\tag 2.7 
$$

Choosing first $U(u) = u^2$ in (2.7), we see that 
$$
\int_\RR u(t)^2 dx + 2 \, \eps \int_0^t \int_\RR  u^2_x dx = \int_\RR u^2_0(x), 
\tag 2.8 
$$
which gives immediately (2.5a) and (2.5b).

Next, we differentiate (2.1) with respect to $x$  and we multiply by $u_x$: 
$$
{1 \over 2} \, \bigl(u^2_x \bigr)_t 
+ u_x \, \bigl(f'(u) \, u_x \bigr)_x 
= 
\eps \, \bigl(u_x \, u_{xx}\bigr)_x - \eps \, u^2_{xx} 
+ \delta \, \bigl(u_x \, u_{xxx} - {1 \over 2} u^2_{xx}\bigr)_x. 
$$ 
Integrating in space, we get 
$$
{1 \over 2} \, {d \over dt} \int_\RR u^2_x \,  dx 
+ \eps \, \int_\RR u^2_{xx} \,  dx 
=     \int_\RR   u_{xx} \, f'(u) \, u_x \, dx 
= - {1 \over 2} \int_\RR f''(u) \, u_x^3 \, dx.
$$
Hence, integrating over some interval $(0,t)$, we find  
$$
\int_\RR  u_x(t)^2 \, dx 
+ 2 \, \eps \, \int_0^t \int_\RR u^2_{xx} \, dx dt
=  \int_\RR  u_{0x}^2 \, dx 
- \int_0^t  \int_\RR f''(u) \, u_x^3 \, dxdt. 
\tag 2.9 
$$

Multiply (2.9) by $\delta$ and add it up with (2.7): 
$$
\aligned 
& \delta \, \int_\RR  u_x(t)^2 \, dx 
+ 2 \, \eps \, \delta \int_0^t \int_\RR u^2_{xx} \, dx dt
\\
& = \int_\RR U(u_0) \, dx - \int_\RR U(u(t)) \, dx 
+ \delta \, \int_\RR  u_{0x}^2 \, dx 
- \eps \, \int_0^t \int_\RR U''(u) \, u^2_x \, dx dt
\\
& \quad - \delta \, \int_0^t  \int_\RR f''(u) \, u_x^3 \, dxdt 
- {\delta \over 2} \, \int_0^t \int_\RR U'''(u) \, u^3_x \, dx dt. 
\endaligned 
$$ 
Choosing $U$ given by 
$$
U(u) = - 2 \int_0^u \bigl(f(v) - f(0)\bigr) \, dv
\tag 2.10 
$$
the last two terms in the above identity cancel out. Since
$$
-c \leq { U''(u) \over 2} \leq c := \|f'\|_\infty, \qquad \text{ for all } u \in \RR, 
$$ 
thus
$$
-c \, u^2 \leq U(u) \leq c \, u^2, \qquad \text{ for all } u \in \RR, 
$$
we finally obtain 
$$ 
\aligned
& \delta \, \int_\RR  u_x(t)^2 \, dx 
+ 2 \, \eps \, \delta \int_0^t \int_\RR u^2_{xx} \, dx dt 
\\
& \leq
\int_\RR c \, u_0^2 \, dx + \int_\RR  c \, u(t)^2 \, dx 
+ \delta \, \int_\RR  u_{0x}^2 \, dx 
+ 2 \, \eps \, \int_0^t \int_\RR c \, u^2_x \, dx dt. 
\endaligned
$$
Hence using (2.8) 
$$
\delta \, \int_\RR  u_x(t)^2 \, dx 
+ 2 \, \eps \, \delta \int_0^t \int_\RR u^2_{xx} \, dx dt 
\leq
2 \, c \, \int_\RR  u_0^2 \, dx + \delta \, \int_\RR  u_{0x}^2 \, dx, 
$$ 
which leads to (2.5c) and (2.5d). The proof of Theorem 2.1 is completed. 
\quad \qed
\enddemo

%______________________________________________________________ 

Recall that by Kruzkov' theory, given $u_0 \in L^2(\RR)$ 
the Cauchy problem (2.3)-(2.4) admits a unique entropy solution 
$u \in L^\infty\bigl(\RR_+, L^2(\RR)\bigr)$' in the sense of Kruzkov's theory. 
See \cite{\refKruzkov, \refDiPerna, \refSzepessy, \refKL}.

\proclaim{Theorem 2.2}
Assume that, for some constant $C_0>0$ independent of $\eps$, 
$$
\|u^\eps_0\|_{L^2(\RR)} + \sqrt{\delta} \, \|u_{0x}^\eps\|_{L^2(\RR)} 
\leq C_0. 
\tag 2.11 
$$ 
\roster 
\item As $\eps \to 0$ with $\delta= O(\eps^2)$ (a subsequence of) 
the solution $u^\eps$ of $(2.1)$-$(2.2)$
converges in $L_{\text{loc}}^p\bigl(\RR_+, L_{\text{loc}}^q(\RR)\bigr)$ 
(for all $1< p < \infty$ and $1 < q < 2$) toward a weak solution of the problem 
$(2.3)$-$(2.4)$.
\item If the stronger condition $\delta= o(\eps^2)$ holds, 
then the limit is the unique entropy solution in the sense of Kruzkov. 
\endroster 
\endproclaim

In Case (1) a subsequence of $u^\eps$ (at least) converges strongly, 
while in Case (2) the whole sequence converges strongly. 
We can conjecture that, in fact, the whole sequence should converge in Case 
(1) as well, but proving it would be very challenging since it  
requires a uniqueness result of nonclassical entropy solutions. 
(See also LeFloch \cite{\refLeFlochtwo}.)

%_________________________________________________________________ 

\demo{Proof} We will apply the general convergence framework established by Schonbek 
\cite{\refSchonbek}. Based on (2.11) and the uniform estimate (2.5a) derived earlier, we can 
select a subsequence of $u^\eps$ converging ``in the sense'' of the Young measures. 
To apply \cite{\refSchonbek}, we only need to control the entropy dissipation measures 
associated with the equation (2.1). Let $U$ be a smooth function with (at most) linear growth at infinity 
and, more precisely, such that $U'$ and $U''$ are uniformly bounded on $\RR$. 
Consider the distribution 
$$
\Gamma^\eps = U(u^\eps)_t + F(u^\eps)_x,  
$$
where as usual $F' = U'\, f'$. With obvious notation consider the decomposition 
$$
\aligned
\Gamma^\eps 
&= \eps \, \bigl(U'(u^\eps) \, u^\eps_x\bigr)_x 
    - \eps \, U''(u^\eps) \, (u^\eps_x)^2 
    + \delta \, \bigl(U'(u^\eps) \, u^\eps_{xx} \bigr)_x 
    - \delta \, U''(u^\eps) \, u^\eps_x \, u_{xx}^\eps
\\
& = \Gamma^\eps_1 + \Gamma^\eps_2 + \Gamma^\eps_3 + \Gamma^\eps_4. 
\endaligned
$$
The estimates below hold for all smooth function $\theta : \RR \times \RR_+ \to \RR$ 
with compact support in $(x,t)$. 

Consider first the term $\Gamma_1^\eps$.  By Cauchy-Schwarz inequality, we get 
$$
\aligned
\bigg|\int_0^\infty \int_\RR \Gamma_1^\eps \, \theta \, dx dt \bigg| 
& = \bigg|\int_0^\infty \int_\RR \eps \, U'(u^\eps) \, u^\eps_x \, \theta_x \, dx dt \bigg| 
\\
& \leq \eps \, C  \, \|u^\eps_x\|_{L^1\bigl(\RR_+, L^2(\RR)\bigr)} 
       \, \|\theta_x\|_{L^\infty\bigl(\RR_+, L^2(\RR)\bigr)} 
\\
& \leq C'  \, \sqrt{\eps} \to 0, 
\endaligned
\tag 2.12i 
$$
where we used (2.5b). This proves that $\Gamma_1^\eps$ converges to zero in the 
sense of distributions. 

Next we simply point out that, by (2.5b) again,
the second term $\Gamma^\eps_2$ remains uniformly bounded in $L^1$: 
$$
\int_0^\infty \int_\RR |\Gamma^\eps_2| \, dx dt 
\leq {1 \over 2} \, \|u^\eps_0\|^2_{L^2(\RR)}.
\tag 2.12ii
$$
To estimate $\Gamma_3$ we use (2.5d):  
$$
\aligned
\bigg|\int_0^\infty \int_\RR \Gamma^\eps_3 \, \theta \, dx dt \bigg| 
& = \bigg|\delta \int_0^\infty \int_\RR U'(u^\eps) \, u^\eps_{xx} \, \theta_x \, dx dt \bigg| 
\\ & \leq
\delta \, C \|u^\eps_{xx}\|_{L^1\bigl(\RR_+, L^2(\RR)\bigr)} \, 
              \|\theta_x\|_{L^\infty\bigl(\RR_+, L^2(\RR)\bigr)} 
\\
& \leq C' \, \sqrt{\delta \over \eps} \to 0,  
\endaligned
\tag 2.12iii 
$$
provided that the mild condition $\delta = o(\eps)$ holds. Therefore $\Gamma^\eps_3$
tends to zero in the sense of distributions.

Finally, we deal with the last term as follows: 
$$
\aligned
\bigg|\int_0^\infty \int_\RR \Gamma^\eps_4 \, \theta \, dx dt \bigg| 
& = \bigg|\int_0^\infty \int_\RR 
\delta \, U''(u^\eps) \, u^\eps_x \, u_{xx}^\eps \, \theta \, dx dt \bigg|\\
& \leq 
\delta \,  C \|u^\eps_{xx}\|_{L^\infty\bigl(\RR_+, L^2(\RR)\bigr)} \, 
              \|u^\eps_x\|_{L^\infty\bigl(\RR_+, L^2(\RR)\bigr)} 
   \,  \|\theta\|_{L^\infty\bigl(\RR\times\RR_+\bigr)}
\\ 
& \leq C' \, {\sqrt{\delta} \over \eps}, 
\endaligned 
\tag 2.12iv 
$$
where we use (2.5b) and (2.5d). 
The upper bound above tends to zero iff $\delta = o(\eps^2)$, in which case we can conclude that 
$\Gamma^\eps_4$ tends to zero in the sense of distributions. 
Under the weaker assumption $\delta = O(\eps^2)$, we see that 
$\Gamma^\eps_4$ is solely bounded in $L^1(\RR \times \RR_+)$ 
as is $\Gamma^\eps_2$.  

The conclusion (1) of the theorem follows immediately from the uniform bounds (2.12) 
by applying Schonbek's convergence theory. Her arguments only show that a subsequence
of $u^\eps$ converges and that the limit is a weak solution of (2.3)-(2.4). 
On the other hand, assuming now the stronger condition  $\delta = o(\eps^2)$
and restricting attention to {\it convex functions\/} $U$, 
in view of (2.12) again and the expression of $\Gamma^\eps_2$ we see 
that the entropy dissipation decomposes in the form 
$$
\Gamma^\eps = \tilde\Gamma^\eps + \Gamma_2^\eps, 
$$ 
where $\tilde \Gamma^\eps \to 0$
in the sense of distributions and $\Gamma^\eps$ is a non-positive bounded measure. 
This shows that all of the entropy inequalities hold in the limit $\eps \to 0$. Thus 
the limit coincides with the unique entropy solution of the problem. 
\quad \qed
\enddemo

%============================================================== 

\heading{3. Multi-Dimensional Conservation Laws} 
\endheading

The estimates and the technique of proof in Section 2 do not apply to multi-dimensional
equations, and markedly different arguments are discussed now. 
Consider the following Cauchy problem: 
$$
u_t + \sum_{j=1}^d f_j(u)_{x_j} 
= \eps \, \sum_{j=1}^d u_{x_j x_j} 
 + \delta \, \sum_{j=1}^d u_{x_j x_j x_j}, \qquad u= u^\eps(x,t), \, x \in \Rd, \, t>0, 
\tag 3.1
$$
$$
u(x,0) = u_0^\eps(x), \quad x \in \Rd. 
\tag 3.2 
$$ 
Provided the initial data $u_0^\eps$ converge weakly to some limit $u_0$ (in $L^2$, say), 
we will now prove that the solutions of (3.1)-(3.2) converge toward the entropy solution of  the 
associated hyperbolic problem: 
$$
u_t + \sum_{j=1}^d f_j(u)_{x_j} = 0, \qquad u= (x,t), \, x \in \Rd, \, t>0, 
\tag 3.3 
$$ 
$$
u(x,0) = u_0(x), \quad x \in \Rd. 
\tag 3.4 
$$ 

Precisely our result are as follows: 

\proclaim{Theorem 3.1} Suppose that the flux-function $f$ is Lipschitz continuous on $\RR$
and that the initial data $u_0$ belong to $L^2(\Rd)$. 
Then the solution $u^\eps$ of $(3.1)$-$(3.2)$ satisfies the following a priori estimates: 
$$
\|u^\eps(t)\|_{L^2(\Rd)} \leq \|u^\eps_0\|_{L^2(\Rd)}, \qquad t \geq 0, 
\tag 3.5a
$$
$$
\sqrt { 2 \, \eps } \, \|u^\eps_x\|_{L^1\bigl(\RR^+,L^2(\Rd)\bigr)}
\leq \|u^\eps_0\|_{L^2(\Rd)}, 
\tag 3.5b 
$$ 
for all $j= 1, \ldots, d$ and $t \geq 0$
$$
\eps \, \|u_{x_j}^\eps(t)\|_{L^2(\Rd)} 
\leq \sqrt{d } \, \|f_j'\|_\infty \, \|u^\eps_0\|_{L^2(\Rd)} 
       + \eps \,  \|u^\eps_{0 x_j}\|_{L^2(\Rd)}
\tag 3.5c
$$
and for all $j,k= 1, \ldots, d$ 
$$
\eps^{3/2} \, \|u^\eps_{x_j x_k}\|_{L^1\bigl(\RR^+,L^2(\Rd)\bigr)}
\leq
 \sqrt{d } \,  \|f_j'\|_\infty \, \|u^\eps_0\|_{L^2(\Rd)} + \eps \,  \|u^\eps_{0 x_j}\|_{L^2(\Rd)}.
\tag 3.5d 
$$
\endproclaim 

For each $u_0 \in L^2(\Rd)$, 
the Cauchy problem (3.3)-(3.4) admits a unique entropy solution 
$u \in L^\infty(\RR_+, L^2(\Rd)$' in the sense of Kruzkov. 
See again \cite{\refKruzkov, \refDiPerna, \refSzepessy, \refKL, \refLeFlochtwo}.

\proclaim{Theorem 3.2}
Assume that, for some constant $C_0>0$ independent of $\eps$, 
$$
\|u^\eps_0\|_{L^2(\Rd)} 
+ \eps \, \sum_{j=1}^d \|u_{0x_j}^\eps\|_{L^2(\RR)} 
\leq C_0. 
\tag 3.6 
$$ 
Then, when $\eps \to 0+$ with $\delta= o(\eps^2)$, the solution $u^\eps$ of $(3.1)$-$(3.2)$
converges in $L_{\text loc}^p\bigl(\RR_+, L_{\text loc}^q(\Rd)\bigr)$ 
(for all $1<p<\infty$ and $1 < q < 2$) toward  
the unique entropy solution in the sense of Kruzkov of the Cauchy problem 
$(3.3)$-$(3.4)$. 
\endproclaim

Recall again that the condition $\delta= o(\eps^2)$ is sharp since, in the opposite case, 
nonclassical solutions violating the Kruzkov entropy inequalities could arise in the limit. 

%---------------------------------------------------------- 

\demo{Proof of Theorem 3.1} We omit the upper-index $\eps$  in the following calculation. 
To derive the $L^2$ bound (3.5a), we multiply the equation (3.1) by $u$ and get 
$$
\aligned 
& \bigl( {|u|^2  \over 2}\bigr)_t  + \sum_{j=1}^d F_j(u)_{x_j}  
\\
& = 
\sum_{j=1}^d \bigl( \eps \, u \, u_{x_j} \bigr)_{x_j} 
- \eps \, \sum_{j=1}^d |u_{x_j}|^2 - {\delta \over 2} \, \sum_{j=1}^d 
\bigl(|u_{x_j }|^2\bigr)_{x_j}
+
\sum_{j=1}^d \bigl( \delta \, u \, u_{x_j x_j} \bigr)_{x_j}, 
\endaligned 
$$ 
where $F_j' = u \, f_j'$ is normalized by the condition $F_j(0)=0$, $j=1, \ldots, d$. 
Integrating over space, we get 
$$ 
{d \over dt} \int_\Rd |u|^2 \, dx = - 2 \, \eps \int_\Rd \sum_{j=1}^d |u_{x_j}|^2 \, dx 
$$ 
and so for all $t \geq 0$ 
$$
\int_\Rd |u(t)|^2 \, dx + 2 \, \eps \, \int_0^t \int_\Rd \sum_{j=1}^d  |u_{x_j}|^2 \, dxdt
= \int_\Rd |u_0|^2 \, dx. 
\tag 3.7 
$$ 

To estimate the gradient of $u$, for $k=1, \cdots, d$ 
we differentiate the equation (3.1) with respect to the variable $x_k$ and 
then multiply by $u_{x_k}$. The right-hand side of (3.1) is linear in $u$ thus 
the calculation for this side is identical to the one we just made,
but with $u$ replaced with $u_{x_k}$. 
On the other hand, the flux term in the left-hand side is nonlinear 
and requires a specific argument: 
$$
{d \over dt} \int_\Rd |u_{x_k}|^2 \, dx 
- \sum_{j=1}^d \int_\Rd 2 \, u_{x_k x_j}\, f_j'(u) \, u_{x_k} \, dx 
=  
- 2 \eps \, \sum_{j=1}^d  \int_\Rd |u_{x_j x_k}|^2 \, dx, 
$$ 
so after integration in time 
$$ 
\aligned 
& \int_\Rd |u_{x_k}(t)|^2 \, dx 
+  
2\, \eps \, \sum_{j=1}^d  \int_0^t\int_\Rd |u_{x_j x_k}|^2 \, dxdt 
\\ 
& \leq 
\int |u_{0x_k}|^2 \, dx 
+ 2 \, \|f_k'\|_\infty \, \sum_{j=1}^d \int_0^t \int_\Rd |u_{x_j x_k}| \, |u_{x_k}| \, dxdt
\\
& \leq \int_\Rd |u_{0x_k}|^2 \, dx 
+ { \|f_k'\|_\infty^2 \over \eps} \, d \, \sum_{j=1}^d \int_0^t \int_\Rd |u_{x_k}|^2 \, dxdt
+ \eps \, \sum_{j=1}^d \int_0^t \int_\Rd |u_{x_j x_k}|^2 \, dxdt. 
\endaligned 
$$ 
Observe that the last term of the right-hand side coincides with the 
last term of the left-hand side. Therefore, 
multiplying the above inequality by $\eps^2$ and using the entropy dissipation 
bound in (3.7), we deduce that  
$$ 
\aligned
& \int_\Rd \eps^2 \, |u_{x_k}(t)|^2 \, dx +  \sum_{j=1}^d 
\int_0^t\int_\Rd \eps^3 \, |u_{x_j x_k}|^2 \, dxdt  
\\
& \leq 
\int_\Rd \eps^2 \, |u_{0 x_k}|^2 \, dx 
+ \|f_k'\|_\infty^2 \, \int_0^t\int_\Rd d \, \eps \, |u_{x_k}|^2 \, dxdt
\\ 
& \leq 
\int_\Rd \eps^2 \, |u_{0 x_k}|^2 \, dx + d \, \|f_k'\|_\infty^2 \, \int_\Rd |u_0|^2 \, dx.  
\endaligned 
\tag 3.8 
$$ 
\quad\qed
\enddemo

%______________________________________________________________ 

\demo{Proof of Theorem 3.2} We will rely on the convergence framework proposed by 
DiPerna \cite{\refDiPerna} 
for $L^\infty$ solutions and generalized to $L^p$ solutions by Szepessy \cite{\refSzepessy}
and Kondo and LeFloch in \cite{\refKL}. 

Consider a Young measure $\nu$ associated with the sequence $u^\eps$ and based on 
the uniform $L^2$ bound (3.5a). (Such Young measures are described in Schonbek \cite{\refSchonbek}). 
To show that $\nu$ is an entropy measure-valued solution, 
we must check entropy inequalities associated with the equation (3.3), 
that is,  
$$ 
\la \nu, U\ra_t + \sum_{j=1}^d \la \nu, F_j \ra_{x_j}  \leq 0, 
\tag 3.9
$$ 
where $U: \RR \to \RR$ is a convex function with (at most )linear growth at infinity 
and the entropy flux $F_j' = U' \, f_j'$ is normalized so that $F_j(0)=0$. 

By the definition of the Young measure, we only need to establish that,
in  the decomposition 
$$ 
\aligned 
& \del_t U(u^\eps) + \sum_{j=1}^d \del_j F_j(u^\eps) 
\\
& = \sum_{j=1}^d 
 \del_j\bigl( \eps \, U'(u^\eps) \, \del_j u^\eps 
+ \delta(\eps) \, U'(u^\eps) \, \del_j^2 u^\eps\bigr) 
\\
& \quad - \sum_{j=1}^d \eps \, U''(u^\eps) \,|\del_j u^\eps|^2 
+ \delta(\eps) \, U''(u^\eps) \, \del_j u^\eps \,  \del_j^2 u^\eps
\\
& =: \Gamma_1^\eps + \Gamma_2^\eps + \Gamma_3^\eps + \Gamma_4^\eps,  
\endaligned 
\tag 4.11 
$$ 
we have 
$$
\Gamma_1^\eps, \Gamma_2^\eps, \Gamma_4^\eps \to 0 
$$
and 
$$
\Gamma_3^\eps \leq 0. 
$$
These convergence properties precisely were established in 
the proof of Theorem 2.2, at least for one-dimensional equations. 
The extension to multi-dimensional equations is immediate in 
view of the uniform estimates (3.5). 
A detailled discussion of the initial condition at $t=0$ 
(which is based on using suitable entropy inequalities) can be found 
in Kondo and LeFloch in \cite{\refKL}. This completes the proof that 
the convergence framework in \cite{\refKL} applies and provides 
the strong convergence toward the unique entropy solution of (3.3)-(3.4). 
\quad\qed
\enddemo 

%============================================================== 

\heading{References} 
\endheading

\item{[\refBLP]} Baiti P., LeFloch P.G., and Piccoli B., 
Uniqueness of classical and nonclassical solutions for nonlinear 
hyperbolic systems, J. Differential Equations 172 (2001), 59--82. 
 
\item{[\refBL]} Bedjaoui N. and LeFloch P.G., 
Diffusive-dispersive traveling waves and kinetic relations I. 
Nonconvex hyperbolic conservation laws, 
J. Differential Equations (2001), to appear.  

\item{[\refCL]} Correia J. and LeFloch P.G., 
Nonlinear diffusive-dispersive limits 
for multidimensional conservation laws, 
in ``Advances in Nonlinear P.D.E.'s and Related Areas", 
A volume in honour of Prof. X. Ding, editors G.Q. Chen et al., 
World Scientific, 1999, pp.~103--123. 

\item{[\refDiPerna]} DiPerna R.J., 
Measure-valued solutions to conservation laws, 
Arch. Rational Mech. Anal. 88 (1985), 223--270. 

\item{[\refHLone]} Hayes B.T. and LeFloch P.G., 
Nonclassical shocks and kinetic relations~: 
scalar conservation laws, 
Arch. Rational Mech. Anal. 139 (1997), 1--56. 

\item{[\refHLtwo]} Hayes B.T. and LeFloch P.G., 
Nonclassical shocks and kinetic relations~: Strictly hyperbolic systems, 
SIAM J. Math. Anal. 31 (2000), 941--991. 

\item{[\refHLthree]} Hayes B.T. and LeFloch P.G., 
Nonclassical shocks and kinetic relations~: Finite difference schemes, 
SIAM J. Numer. Anal. 35 (1998), 2169--2194. 

\item{[\refJMS]} Jacobs D., McKinney W.R. and Shearer M.,
Traveling wave solutions of the modified Korteweg-deVries Burgers equation,
J. Differential Equations 116 (1995), 448--467.

\item{[\refKL]} Kondo C. and LeFloch P.G., 
Measure-valued solutions and well-posedness 
of multi-dimensional conservation laws in a bounded domain, 
Portugal. Math. 58 (2001), 171--194. 

\item{[\refKruzkov]} Kru\v zkov S.N., 
First order quasilinear equations in several independent variables, 
Mat. Sbornik 81 (1970), 285--355; English translation in 
Math. USSR Sb. 10 (1970), 217--243. 

\item{[\refLax]} Lax P.D.,  
The zero dispersion limit, a deterministic analogue of turbulence,
Comm. Pure Appl. Math. 44 (1991), 1047--1056. 

\item{[\refLL]} Lax P.D. and Levermore C.D.,
The small dispersion limit of the Korteweg-deVries equation
Comm. Pure Appl. Math. 36 (1983) I, 253--290, 
II, 571--593, III, 809--829.

\item{[\refLeFlochone]} LeFloch P.G.,
An introduction to nonclassical shocks of systems of conservation laws,  
Proc. International School on Theory and Numerics for Conservation Laws,
Freiburg, Germany, 20-24 Oct. 97, D. Kr\"oner, M. Ohlberger and C. Rohde eds., 
Lectures Notes in Computational Science and Engineering, 
Springer Verlag New York, 1999, pp.~28--73. 

\item{[\refLeFlochtwo]} LeFloch P.G.,
{\it Hyperbolic systems of conservation laws~: 
The theory of classical and nonclassical shock waves,\/}
E.T.H. Lecture Notes Series, Birkh\"auser, 2002.

\item{[\refLN]} LeFloch P.G. and Natalini R., 
Conservation laws with vanishing nonlinear diffusion and dispersion, 
Nonlinear Analysis T.M.A. 36 (1999), 213--230.

\item{[\refLR]} LeFloch P.G and Rohde C., 
High-order schemes, entropy inequalities, and nonclassical shocks,
SIAM J. Numer. Anal. 37 (2000), 2023--2060. 

\item{[\refLT]} LeFloch P.G. and Thanh M.D.,
Nonclassical Riemann solvers and kinetic relations III.
A non-convex hyperbolic model for van der Waals fluids,
Electron. J. Differential Equations 72 (2000), 19 pp. 

\item{[\refMurat]} Murat F., Compacit\'e par compensation, 
Ann. Scuola Norm. Sup. Pisa, Sci. Fis. Mat. 5 (1978), 489--507.

\item{[\refSchonbek]} Schonbek M.E., 
Convergence of solutions to nonlinear dispersive equations, 
Comm. Part. Diff. Eqns. 7 (1982) 959--1000. 

\item{[\refSzepessy]} Szepessy A., 
An existence result for scalar conservation laws using
measure-valued solutions, Comm. Part. Diff. Eqns. 14 (1989) 1329--1350.

\item{[\refTartar]} Tartar L., 
The compensated compactness method applied to systems of conservation laws, 
in ``Systems of Nonlinear Partial Differential Equations'', J.M. Ball ed., 
NATO ASI Series, C. Reidel publishing Col., 1983, pp.~263--285.

\end

\enddocument